\documentclass[leqno, 11pt]{amsart}
\usepackage{amsmath, amsthm, amssymb, latexsym}
  \newlength{\auxwidth}
  \newlength{\auxheight}
\newcounter{theorem}
\newcounter{proposition}
\newcounter{lemma}
\newcounter{corollary}
\newcounter{conjecture}

\newtheorem{definition}{Definition}[section]
\newtheorem{theorem}[definition]{Theorem}
\newtheorem{proposition}[definition]{Proposition}
\newtheorem{lemma}[definition]{Lemma}

\theoremstyle{remark}

\numberwithin{equation}{section}

\newcommand{\op}[1]{\operatorname{#1}}

\newcommand{\brak}[1]{\ensuremath{\langle\! #1\!\rangle}}

\newcommand{\Tr}{\ensuremath{\op{Tr}}}
\newcommand{\tr}{\op{tr}}
\newcommand{\Tra}{\ensuremath{\op{Trace}}}

\newcommand{\Res}{\ensuremath{\op{Res}}}
\newcommand{\res}{\ensuremath{\op{res}}}

\newcommand{\C}{\ensuremath{\mathbb{C}}} 
\newcommand{\bH}{\ensuremath{\mathbb{H}}} 
\newcommand{\N}{\ensuremath{\mathbb{N}}} 
\newcommand{\R}{\ensuremath{\mathbb{R}}} 
\newcommand{\Z}{\ensuremath{\mathbb{Z}}}

\newcommand{\Rd}{\ensuremath{\R^{d+1}}}

\newcommand{\Rdo}{\R^{d+1}\!\setminus\! 0}

\newcommand{\URd}{U\times\R^{d+1}}

\newcommand{\URdo}{U\times(\R^{d+1}\!\setminus\! 0)}

\newcommand{\Ca}[1]{\ensuremath{\mathcal{#1}}}

\newcommand{\cE}{\Ca{E}}
\newcommand{\cF}{\Ca{F}}

\newcommand{\cL}{\ensuremath{\mathcal{L}}}

\newcommand{\cN}{\ensuremath{\mathcal{N}}}

\newcommand{\fg}{\ensuremath{\mathfrak{g}}}
\newcommand{\fh}{\ensuremath{\mathfrak{h}}}

\newcommand{\psivdo}{$\Psi_{H}$DO}
\newcommand{\psivdos}{$\Psi_{H}$DO's}

\newcommand{\pvdo}{\ensuremath{\Psi_{H}}}

\newcommand{\pvdoz}{\ensuremath{\Psi_{H}^{\Z}}}

\newcommand{\psidos}{$\Psi$DO's}

\newcommand{\ord}{{\op{ord}}}

\newcommand{\rk}{\op{rk}}
\newcommand{\im}{\op{im}}

\newcommand{\End}{\ensuremath{\op{End}}}

\renewcommand{\Box}{\square}

\newcommand{\dbarb}{\overline{\partial}_{b}}
\newcommand{\dbarbpq}{\overline{\partial}_{b;p,q}}
\newcommand{\varthetabpq}{\overline{\partial}^{*}_{b;p,q}}

\begin{document}
\title{Szeg\"o projections and new invariants for CR and contact manifolds} 

\author{Rapha\"el Ponge}

\address{Graduate School of Mathematical Sciences, University of Tokyo, Tokyo, Japan.}

\email{ponge@math.ohio-state.edu}

\curraddr{Max Planck Institute for Mathematics, Bonn, Germany.}

\thanks{Research partially supported by NSF grant DMS 0409005}


\begin{abstract}
    In this paper we present the construction in~\cite{Po:AofM1} of several new invariants for CR and contact manifolds by looking at the noncommutative residue traces  
    of various geometric \psivdo\ projections. In the CR setting these operators arise from the $\dbarb$-complex and include the Szeg\"o projections 
    acting on $(p,q)$-forms. 
    In the contact setting they stem from  the generalized Szeg\"o projections at arbitrary integer levels of 
    Epstein-Melrose and from the contact complex of Rumin. In particular, we recover and extend recent results of 
    Hirachi and Boutet de Monvel and answer a question of Fefferman. Furthermore, we give simple algebro-geometric arguments proving 
    that Hirachi's invariant vanishes on strictly pseudoconvex CR manifolds of dimension $4m+1$.
\end{abstract}

 \maketitle 

 \section{Introduction}\label{sec:HI}
Let $D\subset \C^{n+1}$ be a strictly pseudoconvex domain with boundary $\partial D$. Let $\theta$ be a pseudohermitian contact form on $\partial D$, i.e., 
if near a point of $\partial D$ we let 
$\rho(z,\overline{z})$ be a local defining function for $D$ with $\partial \overline{\partial}\rho>0$ then $\theta$ agrees up to a conformal factor 
with $i(\partial -\overline{\partial})\rho$. 

We endow $\partial D$ with the Levi metric defined by the Levi form associated to $\theta$ and we let $S_{\theta}:L^{2}(\partial D)\rightarrow L^{2}(\partial D)$ be the Szeg\"o 
projection associated to this metric and let $k_{S_{\theta}}(z,\overline{w})d\theta^{n}\wedge \theta$ be the Schwartz kernel of $S_{\theta}$. 
As shown by Fefferman~\cite{Fe:BKBMPCD} 
and Boutet de Monvel-Sj\"ostrand~\cite{BS:SSNBS} near the diagonal $w=z$ we can write
\begin{equation}
    k_{S_{\theta}}(z,\overline{w})=\varphi_{\theta}(z,\overline{w})\rho(z,\overline{w})^{-(n+1)}+\psi_{\theta}(z,\overline{w})\log 
    \rho(z,\overline{w}),
     \label{eq:Intro.singularity-Szego}
\end{equation}
 where $\varphi_{\theta}(z,\overline{w})$ and $\psi_{\theta}(z,\overline{w})$ are smooth functions. Then Hirachi defined 
 \begin{equation}
     L(S_{\theta}):=\int_{M} \psi_{\theta}(z,\overline{z})d\theta^{n}\wedge \theta.
 \end{equation}

 \begin{theorem}[Hirachi~\cite{Hi:LSSKGISPD}]
1)  $L(S_{\theta})$ is a CR invariant, i.e., it does not depend on the choice of $\theta$. In particular, this is a biholomorphic invariant of 
$D$.\smallskip

2) $L(S_{\theta})$ is invariant under smooth deformations of the domain $D$.   
 \end{theorem}
 
Subsequently, Boutet de Monvel~\cite{BdM:LTTP} generalized Hirachi's result to the contact setting in terms of the generalized Szeg\"o projections for 
contact manifolds introduced in~\cite{BG:STTO}. 
Such operators are FIO's with complex phase and their kernels admit near the diagonal a singularity similar 
to~(\ref{eq:Intro.singularity-Szego}). In this setting the integral of the leading logarithmic singularity defines a contact invariant. 

It has been asked by Fefferman whether there would exist other invariants like $L(S_{\theta})$, i.e., invariants arising from the integrals of the log 
singularities of geometric operators. The aim of this paper is to explain that there are many such invariants. These invariants can be classified into 
three families:\smallskip

(i) CR invariants coming for the $\dbarb$-complex of Kohn-Rossi~(\cite{KR:EHFBCM}, \cite{Ko:BCM});\smallskip

(ii) Contact invariants arising from the generalized Szeg\"o projections at arbitrary integer level of Epstein-Melrose~\cite{EM:HAITH};\smallskip

(iii) Contact invariants coming from the contact complex of Rumin~\cite{Ru:FDVC}.\smallskip

\noindent The construction of these invariants is based on two main tools:\smallskip

- The Heisenberg calculus of Beals-Greiner~\cite{BG:CHM} and Taylor~\cite{Ta:NCMA};\smallskip

- The noncommutative residue trace for the Heisenberg calculus constructed in~\cite{Po:CRAS1} and~\cite{Po:GAFA1}.\smallskip

To date there is no known example of CR or contact manifold for which one of those invariants is not zero. The only known results are vanishing results: 
Hirachi~\cite{Hi:LSSKGISPD} and Boutet de Monvel~\cite{BdM:LTTP} proved that their 
invariant vanish in dimension $3$, and Boutet de Monvel~\cite{BdM:RIMS05} has announced a proof of the vanishing of the invariant in any dimension, 
but the details of the proof have not appeared yet. In this paper, we mention give simple algebro-geometric arguments proving the vanishing of this
invariant on strictly pseudoconvex CR manifolds of dimension $4m+1$ (see Section~\ref{sec:vanishing}).

The talk is organized as follows. In Section~\ref{sec:HC} we recall few facts about Heisenberg manifolds and the Heisenberg calculus. In Section~\ref{sec:NCR} we 
recall the construction and the main properties of the noncommutative residue for the Heisenberg calculus. In Section~\ref{sec:NCRI} we present the 
construction of the CR invariants for the $\dbarb$-complex. In Section~\ref{sec.Toeplitz} we obtain contact invariants from the generalized Szeg\"o projections 
of Epstein-Melrose. In Section~\ref{sec.contact-complex} we construct contact invariants from Rumin's contact complex. Finally, in Section~\ref{sec:vanishing} 
we establish the vanishing of Hirachi's invariant on strictly pseudoconvex CR manifolds of dimension $4m+1$.
%

\section{Heisenberg calculus}\label{sec:HC}   

 \subsection{Heisenberg manifolds}\label{sec:HM}
A Heisenberg manifold is a pair $(M,H)$ consisting of a manifold $M$ together with a distinguished hyperplane bundle $H 
\subset TM$. Moreover, given another Heisenberg manifold $(M',H')$ we say that a diffeomorphism 
$\phi:M\rightarrow M'$ is a Heisenberg  diffeomorphism when $\phi_{*}H=H'$. 

The main examples of Heisenberg manifolds include the following. 

\emph{a) Heisenberg group}. The $(2n+1)$-dimensional Heisenberg group
$\bH^{2n+1}$ is $\R^{2n+1}=\R \times \R^{2n}$ equipped with the 
group law, 
\begin{equation}
    x.y=(x_{0}+y_{0}+\sum_{1\leq j\leq n}(x_{n+j}y_{j}-x_{j}y_{n+j}),x_{1}+y_{1},\ldots,x_{2n}+y_{2n}).
    \label{eq:Heisenberg.group-law}
\end{equation}
A left-invariant basis for its Lie algebra $\fh^{2n+1}$ is 
provided by the vector-fields, 
\begin{equation}
    X_{0}=\frac{\partial}{\partial x_{0}}, \quad X_{j}=\frac{\partial}{\partial x_{j}}+x_{n+j}\frac{\partial}{\partial 
    x_{0}}, \quad X_{n+j}=\frac{\partial}{\partial x_{n+j}}-x_{j}\frac{\partial}{\partial 
    x_{0}}, 
     \label{eq:Examples.Heisenberg-left-invariant-basis}
\end{equation}
with $j=1,\ldots,n$. For $j,k=1,\ldots,n$ and $k\neq j$ we have the relations,
\begin{equation}
    [X_{j},X_{n+k}]=-2\delta_{jk}X_{0}, \qquad [X_{0},X_{j}]=[X_{j},X_{k}]=[X_{n+j},X_{n+k}]=0.
     \label{eq:Heisenberg.Heisenberg-relations}
\end{equation}
In particular, the subbundle spanned by the vector fields 
$X_{1},\ldots,X_{2n}$ defines a left-invariant Heisenberg structure on 
$\bH^{2n+1}$.\smallskip

 \emph{(b) Codimension 1 foliations.} These are the Heisenberg manifolds $(M,H)$ such that $H$ is integrable in Fr\"obenius' sense, i.e.,
 $C^{\infty}(M,H)$ is closed under the Lie bracket of vector fields.\smallskip

 \emph{(c) Contact manifolds}. A contact manifold is a Heisenberg manifold $(M^{2n+1},H)$ such that near any point of $M$ there exists a contact form 
 annihilating $H$, i.e., a 1-form $\theta$ such that $d\theta_{|_{H}}$ is non-degenerate. 
  When $M$ is orientable it is equivalent to require the 
  existence of a globally defined contact form on $M$  annihilating $H$. More specific examples of contact manifolds include the Heisenberg group 
  $\bH^{2n+1}$, boundaries of strictly pseudoconvex domains $D\subset \C^{2n+1}$, like the sphere $S^{2n+1}$, or even the cosphere bundle 
  $S^{*}M$ of a Riemannian manifold $M^{n+1}$.\smallskip

 \emph{d) Confoliations}. The confoliations of Elyashberg and Thurston in~\cite{ET:C} interpolate between contact manifolds and foliations. They can be seen 
 as oriented Heisenberg manifolds $(M^{2n+1},H)$ together with  a non-vanishing $1$-form $\theta$ on $M$ annihilating $H$ and such that 
$(d\theta)^{n}\wedge \theta\geq 0$.\smallskip

 \emph{e) CR manifolds.} If $D\subset \C^{n+1}$ a bounded domain with boundary $\partial D$ then the maximal complex structure, or CR structure, 
of $T(\partial D)$ is given 
by $T_{1,0}=T(\partial D)\cap T_{1,0}\C^{n+1}$, where $T_{1,0}$ denotes the holomorphic tangent bundle of $\C^{n+1}$. 
More generally, a CR structure on an orientable manifold $M^{2n+1}$ is given by a complex rank $n$ integrable subbundle $T_{1,0}\subset T_{\C}M$ such 
that $T_{1,0}\cap \overline{T_{1,0}}=\{0\}$. 
Besides on boundaries of complex domains, 
such structures naturally appear on real hypersurfaces in $\C^{n+1}$, quotients of the Heisenberg group $\bH^{2n+1}$ by cocompact lattices, 
boundaries of complex hyperbolic spaces, and circle bundles over complex manifolds.  

A real hypersurface $M=\{r=0\}\subset \C^{n+1}$ is  strictly pseudoconvex when the Hessian $\partial \overline{\partial} r$ is positive definite. In general, 
to a CR manifold $M$ we can associate a Levi form $L_{\theta}(Z,W)=-id\theta(Z,\overline{W})$ on the CR tangent bundle $T_{1,0}$ by picking a 
non-vanishing real 1-form $\theta$ annihilating $T_{1,0}\oplus T_{0,1}$.  We then say that $M$ is strictly pseudoconvex 
(resp.~$\kappa$-strictly pseudoconvex) when we can choose $\theta$ so that $L_{\theta}$ is positive definite (resp.~is nondegenerate with $\kappa$ negative 
eigenvalues) at every point. 
 In particular, when this happens  $\theta$ is non-degenerate on $H=\Re (T_{1,0}\oplus 
 T_{0,1})$ and so $(M,H)$ is a contact manifold.

\subsection{Tangent Lie group bundle}
The terminology Heisenberg manifold stems from the fact that the relevant tangent structure in this setting is that of a bundle $GM$ of graded nilpotent Lie 
groups (see~\cite{BG:CHM}, \cite{Be:TSSRG}, \cite{EM:HAITH}, \cite{FS:EDdbarbCAHG},  \cite{Gr:CCSSW}, \cite{Po:Pacific1}, \cite{Ro:INA}). This 
tangent Lie group bundle bundle can be described as follows.

First, we can define an intrinsic Levi form as the 2-form $\cL:H\times H\rightarrow TM/H$ such that, for any point $a 
\in M$ and any sections $X$ and $Y$ of $H$ near $a$, we have 
\begin{equation}
    \cL_{a}(X(a),Y(a))=[X,Y](a) \qquad \bmod H_{a}.
     \label{eq:NCRP.Levi-form}
\end{equation}
In other words the class of $[X,Y](a)$ modulo $H_{a}$ depends only on $X(a)$ and $Y(a)$, not on the germs of $X$ and $Y$ near $a$ (see~\cite{Po:Pacific1}). 

We define the tangent Lie algebra bundle $\fg M$ as the graded Lie algebra bundle consisting of $(TM/H)\oplus H$ together with the 
fields of Lie bracket and dilations such that, for sections $X_{0}$, $Y_{0}$ of $TM/H$ and $X'$, $Y'$ 
of $H$ and for $t\in \R$, we have 
\begin{equation}
    [X_{0}+X',Y_{0}+Y']=\cL(X',Y'), \qquad  t.(X_{0}+X')=t^{2}X_{0}+t X' .
    \label{eq:NCRP.Heisenberg-dilations} 
\end{equation}

Each fiber $\fg_{a}M$, $a \in M$,  is a two-step nilpotent Lie algebra, so by requiring the exponential map to be the identity 
the associated tangent Lie group bundle $GM$ appears as $(TM/H)\oplus H$ together with the grading above and the product law such that, 
for sections $X_{0}$, $Y_{0}$ of $TM/H$ and $X'$, $Y'$ of $H$, we have 
 \begin{equation}
     (X_{0}+X').(Y_{0}+Y')=X_{0}+Y_{0}+\frac{1}{2}\cL(X',Y')+X'+Y'.
 \end{equation}

 Moreover, if $\phi$ is a Heisenberg diffeomorphism from $(M,H)$ onto a Heisenberg manifold $(M',H')$ then, as $\phi_{*}H=H'$ we get linear isomorphisms 
 from $TM/H$ onto $TM'/H'$ and from $H$ onto $H'$ which together give rise to a linear isomorphism 
 $\phi_{H}':TM/H\oplus H\rightarrow TM'/H'\oplus H'$. In fact $\phi_{H}'$ is a graded Lie group isomorphism from $GM$ onto $GM'$ 
 (see~\cite{Po:Pacific1}). 

 On the other hand, we have: 
 
 \begin{proposition}[\cite{Po:Pacific1}]\label{prop:Heisenberg.tangent-Lie-group}
     1) At a point $x \in M$ we have $\rk \cL(x)=2n$ iff $G_{x}M$ is isomorphic to $\bH^{2n+1}\times \R^{d-2n}$.\smallskip
     
     2) If $\dim M=2n+1$ then $(M^{2n+1},H)$ is a contact manifold iff $GM$ is a fiber bundle with typical fiber $\bH^{2n+1}$. 
 \end{proposition}
 
\subsection{Heisenberg calculus}\label{subsec:HC}  
 The Heisenberg calculus is the relevant pseudodifferential calculus to study hypoelliptic 
 operators on Heisenberg manifolds. It was independently introduced by Beals-Greiner~\cite{BG:CHM} and Taylor~\cite{Ta:NCMA} (see 
 also~\cite{BdM:HODCRPDO}, \cite{Dy:POHG}, \cite{Dy:APOHSC}, \cite{EM:HAITH}, \cite{FS:EDdbarbCAHG}, \cite{Po:MAMS1}, 
 \cite{RS:HDONG}).
  
 The initial idea in the Heisenberg calculus, which is due to Stein, is to construct a 
 class of operators on a Heisenberg manifold $(M^{d+1},H)$, called \psivdos, which at any point $a \in M$ are modeled on homogeneous left-invariant 
 convolution operators on the tangent group $G_{a}M$. 
%
 
 Locally the \psidos\ can be described  as follows. Let $U \subset \Rd$ be a local chart together with a frame $X_{0},\ldots,X_{d}$ of $TU$ such that 
 $X_{1},\ldots,X_{d}$ span $H$. Such a chart is called a Heisenberg chart. Moreover, on $\Rd$ we consider the dilations, 
\begin{equation}
     t.\xi=(t^{2}\xi_{0},t\xi_{1},\ldots,t\xi), \qquad \xi\in \Rd, \quad t>0. 
     \label{eq:HC.Heisenberg-dilations}
\end{equation}
 
\begin{definition}1)  $S_{m}(\URd)$, $m\in\C$, is the space of functions 
    $p(x,\xi)$ in $C^{\infty}(U\times\Rdo)$ such that $p(x,t.\xi)=t^m p(x,\xi)$ for any $t>0$.\smallskip

2) $S^m(\URd)$,  $m\in\C$, consists of functions  $p\in C^{\infty}(\URd)$ with
an asymptotic expansion $ p \sim \sum_{j\geq 0} p_{m-j}$, $p_{k}\in S_{k}(\URd)$, in the sense that, for any integer $N$ and 
for any compact $K \subset U$, we have
\begin{equation}
    | \partial^\alpha_{x}\partial^\beta_{\xi}(p-\sum_{j<N}p_{m-j})(x,\xi)| \leq 
    C_{\alpha\beta NK}\|\xi\|^{\Re m-\brak\beta -N}, \qquad  x\in K, \quad \|\xi \| \geq 1,
    \label{eq:NCRP.asymptotic-expansion-symbols}
\end{equation}
where we have let $\brak\beta=2\beta_{0}+\beta_{1}+\ldots+\beta_{d}$ and 
 $\|\xi\|=(\xi_{0}^{2}+\xi_{1}^{4}+\ldots+\xi_{d}^{4})^{1/4}$.
\end{definition}

Next, for $j=0,\ldots,d$ let  $\sigma_{j}(x,\xi)$ denote the symbol (in the 
classical sense) of the vector field $\frac{1}{i}X_{j}$  and set  $\sigma=(\sigma_{0},\ldots,\sigma_{d})$. Then for $p \in S^{m}(\URd)$ we let $p(x,-iX)$ be the 
continuous linear operator from $C^{\infty}_{c}(U)$ to $C^{\infty}(U)$ such that 
    \begin{equation}
          p(x,-iX)f(x)= (2\pi)^{-(d+1)} \int e^{ix.\xi} p(x,\sigma(x,\xi))\hat{f}(\xi)d\xi,
    \qquad f\in C^{\infty}_{c}(U).
    \end{equation}


\begin{definition}
   $\pvdo^{m}(U)$, $m\in \C$, consists of operators $P:C^{\infty}_{c}(U)\rightarrow C^{\infty}(U)$ which are of the form
$P= p(x,-iX)+R$ for some $p$ in $S^{m}(\URd)$, called the symbol of $P$, and some smoothing operator $R$.
\end{definition}

 For any $a\in U$ there is exists a unique affine change of variable $\psi_{a}:\Rd \rightarrow \Rd$ such that $\psi_{a}(a)=0$ 
 and $(\psi_{a})_{*}X_{j}=\frac{\partial}{\partial x_{j}}$ at $x=0$ for $j=0,1,\ldots,d+1$. Then, a continuous operator 
 $P:C^{\infty}_{c}(U)\rightarrow C^{\infty}(U)$ is a \psivdo\ of order $m$ if, and only if, its kernel $k_{P}(x,y)$ has a behavior near the diagonal 
 of the form, 
 \begin{equation}
     k_{P}(x,y) \sim \!\!\sum_{j\geq -(m+d+2)} \!\! (a_{j}(x,\psi_{x}(y))- \!\! \sum_{\brak \alpha =j}\!\! c_{\alpha}(x)\psi_{x}(x)^{\alpha}\log 
     \|\psi_{x}(y)\|), 
      \label{eq:NCRP.kernel-characterization}
 \end{equation}
 with $c_{\alpha}\in C^{\infty}(U)$ and $a_{j}(x,y) \in C^{\infty}(\URdo)$ such that 
 $a_{j}(x,\lambda.y)=\lambda^{j}a_{j}(x,y)$ for any $\lambda>0$. Moreover, $a_{j}(x,y)$ and $c_{\alpha}(x)$, $\brak \alpha=j$, depend only on the symbol 
 of $P$ of degree $-(j+d+2)$.

 
 The class of \psivdos\ is invariant under changes of Heisenberg chart (see~\cite[Sect.~16]{BG:CHM}, \cite[Appendix A]{Po:MAMS1}), so we may 
 extend the definition of \psivdos\ to an arbitrary Heisenberg manifold $(M,H)$ and let them act on sections of a vector bundle $\cE$ over $M$. 
 We let $\pvdo^{m}(M,\cE)$ denote the class of \psivdos\ of order $m$ on $M$ acting on sections 
 of $\cE$. 

Let $\fg^{*}M$ denote the (linear) dual of the Lie algebra bundle $\fg M$ of $GM$ 
 with canonical projection  $\text{pr}: M\rightarrow \fg^{*}M$. As shown 
 in~\cite{Po:MAMS1} (see also~\cite{EM:HAITH}) 
 the principal symbol of an operator $P\in \pvdo^{m}(M,\cE)$ can be intrinsically defined as a symbol $\sigma_{m}(P)$ of the class below. 
 
\begin{definition}
  $S_{m}(\fg^{*}M)$, $m\in \C$, consists of sections $p\in C^{\infty}(\fg^{*}M\setminus 0, \End \textup{pr}^{*}\cE)$ which are homogeneous of 
    degree $m$ with respect to the dilations in~(\ref{eq:NCRP.Heisenberg-dilations}), i.e., we have 
   $p(x,\lambda.\xi)=\lambda^{m}p(x,\xi)$ for any $\lambda>0$. 
\end{definition}

 Next, for any $a \in M$ the convolution 
 on $G_{a}M$ gives rise under the (linear) Fourier transform to a bilinear product for homogeneous symbols, 
 \begin{equation}
     *^{a}: S_{m_{1}}(\fg^{*}_{a}M,\cE_{a})\times S_{m_{2}}(\fg^{*}_{a}M,\cE_{a}) \longrightarrow S_{m_{1}+m_{2}}(\fg^{*}_{a}M,\cE_{a}),
 \end{equation}
 This product depends smoothly on $a$ as much so to yield a product, 
 \begin{gather}
     *: S_{m_{1}}(\fg^{*}M,\cE)\times S_{m_{2}}(\fg^{*}M,\cE) \longrightarrow S_{m_{1}+m_{2}}(\fg^{*}M,\cE),\\
     p_{m_{1}}*p_{m_{2}}(a,\xi)=[p_{m_{1}}(a,.)*^{a}p_{m_{2}}(a,.)](\xi).
     \label{eq:NCRP.product-symbols}
 \end{gather}
 This provides us with the right composition for principal symbols, since we have
 \begin{equation}
     \sigma_{m_{1}+m_{2}}(P_{1}P_{2})=\sigma_{m_{1}}(P_{1})*\sigma_{m_{2}}(P_{2}) \qquad \forall P_{j}\in \pvdo^{m_{j}}(M,\cE).
 \end{equation}
for $P_{1}\in \pvdo^{m_{1}}(M,\cE)$ and $P_{2}\in \pvdo^{m_{2}}(M,\cE)$ such that one of them is properly supported. 

 Notice that when $G_{a}M$ is not commutative, i.e., $\cL_{a}\neq 0$, the product $*^{a}$ is not anymore the pointwise product of symbols and, in 
 particular, is not commutative. 
 Consequently, unless when $H$ is integrable, the product for Heisenberg symbols is not commutative and, while local, it is not microlocal. 
  
 When the principal symbol of $P\in \pvdo^{m}(M,\cE)$ is invertible with respect to the product $*$, the symbolic calculus 
 of~\cite{BG:CHM} allows us to construct a parametrix for $P$ in $\pvdo^{-m}(M,\cE)$. In particular, although not elliptic, 
 $P$ is hypoelliptic with a controlled loss/gain of derivatives (see~\cite{BG:CHM}). 
 
 In general, it may be difficult to determine whether the principal symbol of a given operator $P$ in $\pvdo^{m}(M,\cE)$ is invertible with respect to the product 
 $*$, but  this can be completely determined in terms of a representation theoretic criterion on each tangent 
 group $G_{a}M$, the so-called Rockland condition (see~\cite[Thm.~3.3.19]{Po:MAMS1}).  In particular, 
 if $\sigma_{m}(P)(a,.)$ is \emph{pointwise} invertible with respect to the product $*^{a}$ for any 
 $a \in M$ then  $\sigma_{m}(P)$ is \emph{globally} invertible with respect to $*$. 

\section{Noncommutative residue}\label{sec:NCR}
Let $(M^{d+1}, H)$ be a Heisenberg manifold equipped with a smooth positive density and let $\cE$ be a Hermitian vector bundle over $M$. We 
let $\pvdoz(M,\cE)$ denote the space of \psivdo\ of integer order acting on sections of $\cE$. 

\subsection{Logarithmic singularity} 
 Let $P:C^{\infty}(M,\cE)\rightarrow C^{\infty}(M,\cE)$ be a \psivdo\ of integer order $m$. Then it follows from~(\ref{eq:NCRP.kernel-characterization}) 
 that in a trivializing Heisenberg chart the kernel $k_{P}(x,y)$ of $P$ has a behavior near the diagonal of the form,  
      \begin{equation}
        k_{P}(x,y)=\!\!\sum_{-(m+d+2)\leq j\leq 1}\!\!a_{j}(x,-\psi_{x}(y)) - c_{P}(x)\log \|\psi_{x}(y)\| + 
        \op{O}(1),
        \label{eq:NCRHC.log-singularity}
    \end{equation}
    where $a_{j}(x,y)$ is homogeneous of degree $j$ in $y$ with respect to the dilations~(\ref{eq:HC.Heisenberg-dilations}). Furthermore, we have
    \begin{equation}
        c_{P}(x)=|\psi_{x}'|\int_{\|\xi\|=1}p_{-(d+2)}(x,\xi)d\xi,
         \label{eq:NCRHC.formula-cP}
    \end{equation}
    where $p_{-(d+2)}(x,\xi)$ is the homogeneous symbol of degree $-(d+2)$ of $P$. 
    
    Let $|\Lambda|(M)$ be the bundle of densities on $M$. Then we have:
\begin{proposition}[\cite{Po:CRAS1}, \cite{Po:GAFA1}]
 The coefficient $c_{P}(x)$ makes sense intrinsically on $M$ as a section of $|\Lambda|(M)\otimes \End \cE$.
\end{proposition}

\subsection{Noncommutative residue}
From now on we assume $M$ compact. Therefore, for any $P\in \pvdoz(M,\cE)$ we can let
    \begin{equation}
        \Res P= \int_{M} \tr_{\cE} c_{P}(x).
         \label{eq:NCR}
    \end{equation}

 If $P$ is in $\pvdo^{m}(M,\cE)$ 
with $\Re m<-(d+2)$ then $P$ is trace-class. It can be shown that we have an analytic continuation of the 
trace to \psivdos\ of non-integer orders which is analogous to that for classical \psidos\ in~\cite{KV:GDEO}. Moreover, on \psivdos\ of integer 
orders this analytic extension of the trace induces a residual functional agreeing with~(\ref{eq:NCR}), so that we have: 

\begin{proposition}
    Let $P\in \pvdoz(M,\cE)$. Then for any family $(P(z))_{z \in \C}\subset \pvdo^{*}(M,\cE)$ which is holomorphic in the sense of~\cite{Po:MAMS1} 
    and such that $P(0)=P$ and $\ord P(z)=z+\ord P$ we have
    \begin{equation}
        \Res P=-\res_{z=0} \Tra P(z).
    \end{equation}
\end{proposition}

Thus the functional~(\ref{eq:NCR}) is the analogue for the Heisenberg calculus of the noncommutative residue of 
Wodzicki~(\cite{Wo:LISA}, \cite{Wo:NCRF}) and Guillemin~\cite{Gu:NPWF}. Furthermore, we have:

\begin{proposition}[\cite{Po:CRAS1}, \cite{Po:GAFA1}]
1) Let $\phi$ be a Heisenberg diffeomorphism from $(M,H)$ onto a Heisenberg manifold $(M',H')$. Then for any $P \in \pvdoz(M,\cE)$ we have $\Res 
     \phi_{*}P=\Res P$.\smallskip  

2) $\Res$ is a trace on the algebra $\pvdoz(M,\cE)$ which vanishes on differential operators and on \psivdos\ of integer order~$\leq -(d+3)$.\smallskip  
  
3) If $M$ is connected then $\Res$ is the unique trace up to constant multiple.
\end{proposition}

Let $D\subset \C^{n+1}$ be a strictly pseudoconvex domain with boundary $\partial D$ and let $\theta$ be a pseudohermitian contact form on $\partial D$. 
We endow $\partial D$ with the associated Levi metric and we let $S_{\theta}:L^{2}(\partial D)\rightarrow L^{2}(\partial D)$ be the corresponding 
Szeg\"o projection. Then $S_{\theta}$ is a \psivdo\ of order 0 and with the notation of~(\ref{eq:Intro.singularity-Szego})  
we have $c_{S_{\theta}}(z)=-\frac{1}{2}\psi_{\theta}(z,\overline{z})d\theta^{n}\wedge \theta$. Thus, 
\begin{equation}
    \Res S_{\theta}=-\frac{1}{2}L(S_{\theta}).
\end{equation}
This shows that Hirachi's invariant can be interpreted as a noncommutative residue. 

 \section{CR invariants from the $\dbarb$-complex}\label{sec:NCRI}
Let $M^{2n+1}$ be a compact orientable CR manifold with CR tangent bundle $T_{1,0}\subset T_{\C}M$, so that $H=\Re (T_{1,0}\oplus T_{0,1})\subset TM$ is 
a hyperplane bundle of $TM$ admitting an (integrable) complex structure. Let $\theta$ be a global non-zero real 1-form  annihilating $H$ and let 
$L_{\theta}$ be the associated Levi form,  
 \begin{equation}
    L_{\theta}(Z,W)=-id\theta(Z,\overline{W})=i\theta([Z,\overline{W}]), \qquad Z,W\in C^{\infty}(M,T_{1,0}). 
      \label{eq:CR.Levi-form}
 \end{equation}

Let $\cN$ be a supplement of $H$ in $TM$. 
This is an orientable line bundle which gives rise to the splitting, 
\begin{equation}
    T_{\C}M=T_{1,0}\oplus T_{0,1}\oplus (\cN\otimes \C).
    \label{eq:CR-decomposition}
\end{equation}
Let $\Lambda^{1,0}$ and $\Lambda^{0,1}$ denote the annihilators 
in $T^{*}_{\C}M$ of $T_{0,1}\oplus (\cN\otimes \C)$ and $T_{1,0}\oplus (\cN\otimes \C)$ respectively and 
for $p,q=0,\ldots,n$ let $\Lambda^{p,q}=(\Lambda^{1,0})^{p}\wedge (\Lambda^{0,1})^{q}$ be the bundle of $(p,q)$-forms. 
Then we have the splitting, 
\begin{equation}
    \Lambda^{*}T_{\C}^{*}M=(\bigoplus_{p,q=0}^{n}\Lambda^{p,q})\oplus \theta\wedge  \Lambda^{*}T_{\C}^{*}M.
     \label{eq:CR-Lambda-pq-decomposition}
\end{equation}
Notice that this decomposition does not depend on the choice of $\theta$, but it does depend on that of $\cN$. 

The complex $\dbarb:C^{\infty}(M, \Lambda^{p,*})\rightarrow C^{\infty}(M,\Lambda^{p,*+1})$ of 
Kohn-Rossi~(\cite{KR:EHFBCM},~\cite{Ko:BCM}) is defined as follows. For any $\eta \in C^{\infty}(M, \Lambda^{p,q})$ we can uniquely decompose $d\eta$ as 
\begin{equation}
    d\eta =\dbarbpq \eta + \partial_{b;p,q}\eta + \theta \wedge \cL_{X_{0}}\eta,
     \label{eq:CR.dbarb}
\end{equation}
where $\dbarbpq \eta $  and $\partial_{b;p,q}\eta$ are sections of $\Lambda^{p,q+1}$ and $\Lambda^{p+1,q}$ respectively and $X_{0}$ is the section of $\cN$ 
such that $\theta(X_{0})=1$. Thanks to the integrability of $T_{1,0}$ we have $\overline{\partial}_{b;p,q+1}\circ \dbarbpq=0$, so that we get a chain 
complex. Notice that this complex depends on the CR structure of $M$ and on the choice of $\cN$. 

Assume now that $M$ is endowed with a Hermitian metric $h$ on $T_{\C}M$ which commutes with complex conjugation and 
makes the splitting~(\ref{eq:CR-decomposition}) become orthogonal. The associated Kohn Laplacian is
\begin{equation}
    \Box_{b;p,q}= \overline{\partial}^{*}_{b;p,q+1}\dbarbpq +\overline{\partial}_{b;p,q-1}\overline{\partial}_{b;p,q}^{*}.
\end{equation}

For $x \in M$ let $\kappa_{+}(x)$ (resp.~$\kappa_{-}(x)$) be the number of positive (resp.~negative) eigenvalues of $L_{\theta}$ at $x$. We then say 
that the condition $Y(q)$ holds when at every point $x \in M$ we have
\begin{equation}
    q \not \in \{\kappa_{-}(x),\ldots,n-\kappa_{+}(x)\}\cup  \{\kappa_{+}(x),\ldots,n-\kappa_{-}(x)\}.
\end{equation}
For instance, when $M$ is $\kappa$-strictly pseudoconvex we have $\kappa_{-}(x)=n-\kappa_{+}(x)=\kappa$, so the condition $Y(q)$ exactly means that we 
must have $q \neq \kappa$ and $q\neq n-\kappa$. 

\begin{proposition}[see~{\cite[Sect.~21]{BG:CHM}}, {\cite[Sect.~3.5]{Po:MAMS1}}]\label{prop:CR.Kohn-Laplacian}
    The Kohn Laplacian $\Box_{b;p,q}$ admits a parametrix in $\pvdo^{-2}(M,\Lambda^{p,q})$ iff the condition $Y(q)$ is satisfied. 
\end{proposition}

Let $S_{b;p,q}$ be the Szeg\"o projection on $(p,q)$-forms, i.e., the orthogonal projection onto $\ker \Box_{b;p,q}$. We also consider 
the orthogonal projections  $\Pi_{0}(\dbarbpq)$ and $\Pi_{0}(\varthetabpq)$ onto $\ker 
\dbarbpq$ and $\ker \varthetabpq=(\im \overline{\partial}_{b;p,q-1})^{\perp}$. In fact, as $\ker \dbarbpq=\ker \Box_{b;p,q}\oplus \im 
\overline{\partial}_{b;p,q-1}$ we have $\Pi_{0}(\dbarbpq)=S_{b;p,q}+1-\Pi_{0}(\varthetabpq)$, that is, 
\begin{equation}
    S_{b;p,q}= \Pi_{0}(\dbarbpq)+\Pi_{0}(\varthetabpq)-1.
     \label{eq:Szego-projections-dbarb-varthetab}
\end{equation}

Let $N_{b;p,q}$ be the partial inverse of $\Box_{b;p,q}$, so that $N_{b;p,q}\Box_{b;p,q}=\Box_{b;p,q}N_{b;p,q}=1-S_{b;p,q}$. Then it can be shown
(see, e.g.,~\cite[pp.~170--172]{BG:CHM}) that we have
\begin{equation}
    \Pi_{0}(\dbarbpq)=1-\overline{\partial}^{*}_{b;p,q}N_{b;p,q+1}\dbarbpq, \qquad  
    \Pi_{0}(\varthetabpq)=1-\overline{\partial}_{b;p,q-1}N_{b;p,q-1} \overline{\partial}_{b;p,q-1}^{*}.
     \label{eq:projections-dbarb-varthetab}
\end{equation}

By Proposition~\ref{prop:CR.Kohn-Laplacian} when the condition $Y(q)$ holds at every point the operator $\Box_{b;p,q}$  
admits a parametrix in $\pvdo^{-2}(M,\Lambda^{p,q})$ and then $S_{b;p,q}$ is a smoothing operator and $N_{b;p,q}$ is a \psivdo\ of 
order $-2$. Therefore, using~(\ref{eq:projections-dbarb-varthetab}) 
we see  that if the condition $Y(q+1)$ (resp.~$Y(q-1)$) holds 
everywhere then $\Pi_{0}(\dbarbpq)$ (resp.~$\Pi_{0}(\varthetabpq)$) is a \psivdo. 

Furthermore, in view of~(\ref{eq:Szego-projections-dbarb-varthetab}) we also see
 that if at every point the condition $Y(q)$ fails, but the conditions $Y(q-1)$ and 
$Y(q+1)$ hold, then the Szeg\"o projection $S_{b;p,q}$ is a zero'th order \psivdo\ projection. Notice that this may happen if, and only if, $M$ is 
$\kappa$-strictly pseudoconvex with $\kappa=q$ or $\kappa=n-q$.
 
Bearing all this in mind we have: 

\begin{theorem}[\cite{Po:AofM1}]\label{CR:Thm1}
1)  The following noncommutative residues are CR diffeomorphism invariants of $M$: \smallskip
    
    (i) $\Res \Pi_{0}(\dbarbpq)$ when the condition $Y(q+1)$ holds everywhere;\smallskip 
    
    (ii) $\Res \Pi_{0}(\varthetabpq)$ when the condition $Y(q-1)$ holds everywhere;\smallskip
    
    (iii) $\Res S_{b;p,\kappa}$ and $\Res S_{b;p,n-\kappa}$ when $M$ is $\kappa$-strictly pseudoconvex.\smallskip
    
\noindent In particular, they depend neither on the choice of the line bundle 
$\cN$, nor on that of the Hermitian metric $h$.\smallskip

2) The noncommutative residues (i)--(iii) are invariant under deformations of the CR structure coming from deformations of the complex structure of 
$H$. 
\end{theorem}

Specializing Theorem~\ref{CR:Thm1} to the strictly pseudoconvex case we get: 

\begin{theorem}[\cite{Po:AofM1}]
    Suppose that $M$ is a compact strictly pseudoconvex CR manifold. Then:\smallskip
    
    1) $\Res S_{b;p,j}$, $j=0,n$, and $\Res \Pi_{0}(\dbarbpq)$, $q=1,\ldots,n-1$,  
are CR diffeomorphism invariants of $M$.  In particular, when $M$ is the boundary of a strictly pseudoconvex domain 
$D\subset \C^{n}$ they give rise to bilholomorphism invariants of $D$.\smallskip    

2) The above residues are invariant under deformations of the CR structure. 
\end{theorem}

\section{Invariants of generalized Szeg\"o projections}\label{sec.Toeplitz}
Let $(M^{2n+1},H)$ be an orientable contact manifold.  
Given a contact form $\theta$ on $M$ annihilating $H$ we let $X_{0}$ be the Reeb 
vector field of $\theta$, i.e., the unique vector field $X_{0}$ such that $\iota_{X_{0}}\theta=1$ and $\iota_{X_{0}}d\theta=0$. 
In addition, we let $J$ be an 
 almost complex structure on $H$  which is calibrated in the sense that $d\theta(X,JX)>0$ for any nonzero section $X$ of $H$. Extending $J$ to $TM$ by 
 requiring to have  
$JX_{0}=0$, we can equip $TM$ with the Riemannian metric $g_{\theta, J}=d\theta(.,J.)+\theta^{2}$.

In this context Szeg\"o projections have been defined by Boutet de Monvel and Guillemin in~\cite{BG:STTO} as an FIO with complex phase. 
This construction has been further generalized by Epstein-Melrose~\cite{EM:HAITH}  as follows. 

Let $\bH^{2n+1}$ be the Heisenberg group of dimension $2n+1$ consisting of $\R^{2n+1}$ together with the 
group law~(\ref{eq:Heisenberg.group-law}). 
Let $\theta^{0}=dx_{0}+\frac{1}{2}\sum_{j=1}^{n}(x_{j}dx_{n+j}-x_{n+j}dx_{j})$ be the standard left-invariant contact form of 
$\bH^{2n+1}$; its Reeb 
vector field is $X^{0}_{0}=\frac{\partial}{\partial x_{0}}$. 

 For $j=1,\ldots,n$  let $X^{0}_{j}=\frac{\partial}{\partial x_{j}}+\frac{1}{2}x_{n+j}\frac{\partial}{\partial 
    x_{0}}$ and $X^{0}_{n+j}=\frac{\partial}{\partial x_{n+j}}-\frac{1}{2}x_{j}\frac{\partial}{\partial 
    x_{0}}$ then $X_{1}^{0},\ldots,X_{2n}^{0}$ form a left-invariant frame of $H^{0}=\ker \theta^{0}$ and sastify the 
    relations~(\ref{eq:Heisenberg.Heisenberg-relations}). The standard CR structure of $\bH^{2n+1}$ is then given by the complex structure $J^{0}$ on $H^{0}$ 
    such that $J^{0}X_{j}^{0}=X_{n+j}^{0}$ and 
 $J^{0}X_{n+j}=-X_{j}$.  Moreover, it follows from~(\ref{eq:Heisenberg.Heisenberg-relations}) 
 that $J^{0}$ is calibrated  with respect to $\theta^{0}$ and that 
 $X_{0}^{0}, X_{1}^{0},\ldots,X_{2n}^{0}$ form an orthonormal frame of 
 $T\bH^{2n+1}$ with respect to the metric $g_{\theta^{0},J^{0}}$. 
 
 The scalar Kohn Laplacian  on $\bH^{2n+1}$ is equal to
\begin{equation}
     \Box_{b,0}^{0}=-\frac{1}{2}((X_{1}^{0})^{2}+\ldots+(X_{2n}^{0})^{2})+i\frac{n}{2}X_{0}^{0}.
\end{equation}
For $\lambda \in \C$ the operator $-\frac{1}{2}((X^{0}_{1})^{2}+\ldots+(X^{0}_{2n})^{2})+i \lambda X^{0}_{0}$ is invertible if, and only if,  we have $\lambda \not \in  
\pm(\frac{n}{2}+\N)$ (see~\cite{FS:EDdbarbCAHG}, \cite{BG:CHM}). For $k=0,1,\ldots$ the orthogonal projection 
$\Pi_{0}(\Box_{b}+ikX^{0}_{0})$ onto the kernel of $\Box_{b}+ikX^{0}_{0}$ is a left-invariant homogeneous  \psivdo\ of order~$0$ 
(see~\cite[Thm.~6.61]{BG:CHM}).  We then let $s_{k}^{0}\in S_{0}((\fh^{2n+1})^{*})$ 
denote its symbol, so that  we have $\Pi_{0}(\Box_{b}+ikX^{0}_{0})=s_{k}^{0}(-iX^{0})$. 

Next, since $(M,H)$ is a contact manifold by Proposition~\ref{prop:Heisenberg.tangent-Lie-group} 
the tangent Lie group bundle $GM$ is a fiber bundle with typical fiber 
$\bH^{2n+1}$. A local trivialization near a given point $a \in M$ is obtained as follows. 

Let $X_{1},\ldots,X_{2n}$ be a local orthonormal frame of $H$ on an open neighborhood $U$ of $a$ and which is admissible in the sense 
that $X_{n+j}=JX_{j}$ for $j=1,\ldots,n$. In addition, let $\underline{X_{0}(a)}$ denote the 
class of $X_{0}(a)$ in $T_{a}M/H_{a}$. Then as shown in~\cite{Po:Pacific1} 
the map $\phi_{X,a}:(T_{a}M/H_{a})\oplus H_{a}\rightarrow \R^{2n+1}$ such that 
\begin{equation}
    \phi_{X,a}(x_{0}\underline{X_{0}(a)}+x_{1}X_{1}(a)+\ldots + 
    x_{2n}X_{2n}(a))=(x_{0},\ldots,x_{2n}), \qquad x_{j}\in \R, 
\end{equation}
gives rise to a Lie group isomorphism from $G_{a}M$ onto $\bH^{2n+1}$. In fact, as $\phi_{X,a}$ depends smoothly on $a$ we get a  
fiber bundle trivialization of $GM|_{U}\simeq U\times \bH^{2n+1}$. 

For $j=0,\ldots,2n$ let $X_{j}^{a}$ be the model vector field of $X_{j}$ at $a$ as defined in~\cite{Po:Pacific1}. This is 
the unique left-invariant vector field on $G_{a}M$ which, in the coordinates provided by $\phi_{X,a}$, agrees with $\frac{\partial}{\partial x_{j}}$ at 
$x=0$. Therefore, we have $X_{j}^{a}=\phi_{X,a}^{*}X_{j}^{0}$ and so we get
\(
    \phi_{X,a}^{*}\Box_{b}^{0}=-\frac{1}{2}((X_{1}^{a})^{2}+\ldots+(X_{2n}^{a})^{2})+i\frac{n}{2}X_{0}^{a}. 
\)

If $\tilde{X}_{1},\ldots,\tilde{X}_{2n}$ is another admissible orthonormal frame of $H$ near $a$, then we pass from 
$(\tilde{X}_{1}^{a},\ldots,\tilde{X}_{2n}^{a})$ to $(X_{1}^{a},\ldots,X_{2n}^{a})$ by an orthogonal linear transformation, which leaves the 
expression $(X_{1}^{a})^{2}+\ldots+(X_{2n}^{a})^{2}$ unchanged. Therefore, the differential operator $\Box_{b}^{a}:=\phi_{X,a}^{*}\Box_{b}^{0}$ makes sense 
independently of the choice of the admissible frame $X_{1},\ldots,X_{2n}$ near $a$. 

On the other hand, as $\phi_{X,a}$ induces a unitary transformation from $L^{2}(G_{a}M)$ onto $L^{2}(\bH^{2n+1})$ we have 
$\Pi_{0}(\Box_{b}^{a}+ikX_{0}^{a})=\Pi_{0}(\phi_{X,a}^{*}(\Box_{b}^{0}+ikX_{0}^{0}))=\phi_{X,a}^{*}\Pi_{0}(\Box_{b}^{0}+ikX_{0}^{0}))$. 
 Hence $\Pi_{0}(\Box_{b}^{a}+ikX_{0}^{a})$ is a zero'th order left-invariant homogeneous \psivdo\ on $G_{a}M$ with symbol 
$s_{k}^{a}(\xi)=\phi_{X,a}^{*}s_{k}^{0}(\xi)=s_{k}^{0}((\phi_{X,a}^{-1})^{t}\xi)$. 
 In fact, since $\phi_{X,a}$ depends smoothly on $a$ we obtain:

\begin{proposition}
    For $k=0,1,\ldots$ there is a uniquely defined symbol $s_{k}\in S_{0}(\fg^{*}M)$ such that, for any admissible orthonormal frame $X_{1},\ldots,X_{d}$ of $H$ 
    near a 
    point $a \in M$, we have $s_{k}(a,\xi)=\phi_{X,a}^{*}s_{k}^{0}(\xi)$ for any $(a,\xi)\in \fg^{*}M\setminus 0$.
\end{proposition}

We call $s_{k}$ the \emph{Szeg\"o symbol at level} $k$. This definition \emph{a priori} depends on the contact form $\theta$ and the almost complex 
structure $J$, but we have: 

\begin{lemma}[\cite{EM:HAITH}, \cite{Po:AofM1}]
    (i) The symbol $s_{k}$ is invariant under conformal changes of contact form.\smallskip
    
    (ii) The change $(\theta, J)\rightarrow (-\theta, -J)$ transforms $s_{k}$ into $s_{k}(x,-\xi)$.\smallskip
    
    (iii) The symbol $s_{k}$ depends on $J$ only up to homotopy of idempotents in $S_{0}(\fg^{*}M)$.
 \end{lemma}
 
From now on we let $\cE$ be a Hermitian vector bundle over $M$. 
\begin{definition}[{\cite[Chap.~6]{EM:HAITH}}]\label{def:IGSP-definition}
    For $k=0,1,\ldots$ a generalized Szeg\"o projection at level $k$ is a \psivdo\ projection $S_{k}\in \pvdo^{0}(M,\cE)$ with principal symbol 
    $s_{k}\otimes \op{id}_{\cE}$.
\end{definition}

Generalized Szeg\"o projections at level $k$ always exist (see~\cite{EM:HAITH}, \cite{Po:AofM1}). Moreover, when 
$k=0$  and $\cE$ is the trivial line bundle the above definition allows us to recover the Szeg\"o projections of~\cite{BG:STTO} 
(see~\cite{Po:AofM1}).  In particular, when $M$ is strictly pseudoconvex the Szeg\"o projection $S_{b,0}$ is a generalized Szeg\"o projection at level $0$. 

Given a generalized Szeg\"o projection at level $k$ we define 
\begin{equation}
    L_{k}(\cE)=\Res S_{k}.
\end{equation}
In fact, we have: 

\begin{proposition}[\cite{Po:AofM1}]
The value of  $L_{k}(\cE)$ does not depend on the choice of $S_{k}$.
\end{proposition}

Next, recall that the $K$-group $K^{0}(M)$ can be described as the group of formal differences of stable homotopy classes of (smooth) vector bundles over $M$, 
where a stable homotopy between vector bundles $\cE_{1}$ and $\cE_{2}$ is given by an auxiliary vector bundle $\cF$ and a vector bundle isomorphism 
$\phi:\cE_{1}\oplus \cF \simeq \cE_{2}\oplus \cF$. Then we obtain:

\begin{theorem}[\cite{Po:AofM1}]
1) $L_{k}(\cE)$ depends only on the Heisenberg diffeomorphism class of $M$ and on the $K$-theory class of $\cE$. In particular, it depends neither 
    on the contact form $\theta$, nor on the almost complex structure $J$.\smallskip 
    
2) $L_{k}(\cE)$ invariant is under deformations of the contact structure. 
\end{theorem}


\section{Invariants from the contact complex}\label{sec.contact-complex}
Let $(M^{2n+1},H)$ be an orientable contact manifold. Let $\theta$ be a contact form on $M$ and let $X_{0}$ be its 
Reeb vector field of $\theta$. We also let $J$ be a calibrated almost complex structure on $H$ and  we endow $TM$ with the Riemannian metric 
$g_{\theta,J}=d\theta(.,J.)+\theta^{2}$.

Observe that the splitting $TM=H\oplus \R X_{0}$ allows us to identify  
$H^{*}$ with the annihilator of $X_{0}$ in $T^{*}M$. More generally, identifying $\Lambda^{k}_{\C}H^{*}$ with $\ker \iota_{X_{0}}$, where 
$\iota_{X_{0}}$ denotes the contraction operator by $X_{0}$, gives the splitting
\begin{equation}
    \Lambda^{*}_{\C}TM=(\bigoplus_{k=0}^{2n}\Lambda^{k}_{\C}H^{*}) \oplus (\bigoplus_{k=0}^{2n} \theta\wedge \Lambda^{k}_{\C}H^{*}).
     \label{eq:contact.decomposition-forms}
\end{equation}
 
For any horizontal form $\eta\in C^{\infty}(M,\Lambda^{k}_{\C}H^{*})$ we can write
$d\eta= d_{b}\eta+\theta \wedge \cL_{X_{0}}\eta$,
where $d_{b}\eta$ is the component of $d\eta$ in $\Lambda^{k}_{\C}H^{*}$. This does not provide us with a complex, for we have 
$d_{b}^{2}=-\cL_{X_{0}}\varepsilon(d\theta)=-\varepsilon (d\theta)\cL_{X_{0}}$, where $\varepsilon(d\theta)$ denotes the exterior multiplication 
by $d\theta$.

The contact complex of Rumin~\cite{Ru:FDVC} is an attempt to get 
a complex of horizontal differential forms by forcing the equalities $d_{b}^{2}=0$ and $(d^{*}_{b})^{2}=0$.

A natural way to modify $d_{b}$ to get the equality $d_{b}^{2}=0$ is to restrict 
$d_{b}$ to the subbundle $\Lambda^{*}_{2}:=\ker \varepsilon(d\theta) \cap \Lambda^{*}_{\C}H^{*}$, since the latter
is closed under $d_{b}$ and is annihilated by $d_{b}^{2}$. 

Similarly, we get the equality $(d_{b}^{*})^{2}=0$ by restricting $d^{*}_{b}$ to the subbundle 
$\Lambda^{*}_{1}:=\ker \iota(d\theta)\cap \Lambda^{*}_{\C}H^{*}=(\im \varepsilon(d\theta))^{\perp}\cap \Lambda^{*}_{\C}H^{*}$, where 
$\iota(d\theta)$ denotes the interior product 
with $d\theta$. This amounts to replace $d_{b}$ by $\pi_{1}\circ d_{b}$, where $\pi_{1}$ is the orthogonal projection onto $\Lambda^{*}_{1}$.

In fact, since $d\theta$ is nondegenerate on $H$ the operator $\varepsilon(d\theta):\Lambda^{k}_{\C}H^{*}\rightarrow \Lambda^{k+2}_{\C}H^{*}$  is 
injective for $k\leq n-1$ and surjective for $k\geq n+1$. This implies that $\Lambda_{2}^{k}=0$ for $k\leq n$ and $\Lambda_{1}^{k}=0$ for $k\geq n+1$. 
Therefore, we only have two halves of complexes. 

As observed by Rumin~\cite{Ru:FDVC} we get a full complex by connecting 
the two halves by means of the operator $D_{R,n}:C^{\infty}(M,\Lambda_{\C}^{n}H^{*}) \rightarrow 
C^{\infty}(M,\Lambda_{\C}^{n}H^{*})$ such that 
\begin{equation}
    D_{R,n}=\cL_{X_{0}}+d_{b,n-1}\varepsilon(d\theta)^{-1}d_{b,n},
\end{equation}
where $\varepsilon(d\theta)^{-1}$ is the inverse of $\varepsilon(d\theta):\Lambda^{n-1}_{\C}H^{*}\rightarrow \Lambda^{n+1}_{\C}H^{*}$. Notice that 
$D_{R,n}$ is a second order differential 
operator.  This allows us to get the contact complex, 
\begin{equation}
    C^{\infty}(M)\stackrel{d_{R,0}}{\rightarrow}
    \ldots 
    C^{\infty}(M,\Lambda^{n})\stackrel{D_{R,n}}{\rightarrow} C^{\infty}(M,\Lambda^{n}) 
    \ldots \stackrel{d_{R,2n-1}}{\rightarrow} C^{\infty}(M,\Lambda^{2n}).
     \label{eq:contact-complex}
\end{equation}
where $d_{R,k}$ agrees with $\pi_{1}\circ d_{b}$ for $k=0,\ldots,n-1$ and  with $d_{R,k}=d_{b}$ otherwise. 

The contact Laplacian is defined as follows. In degree $k\neq n$ this is the differential operator 
$\Delta_{R,k}:C^{\infty}(M,\Lambda^{k})\rightarrow C^{\infty}(M,\Lambda^{k})$ such that
\begin{equation}
    \Delta_{R,k}=\left\{
    \begin{array}{ll}
        (n-k)d_{R,k-1}d^{*}_{R,k}+(n-k+1) d^{*}_{R,k+1}d_{R,k},& \text{$k=0,\ldots,n-1$},\\
         (k-n-1)d_{R,k-1}d^{*}_{R,k}+(k-n) d^{*}_{R,k+1}d_{R,k},& \text{$k=n+1,\ldots,2n$}.
         \label{eq:contact-Laplacian1}
    \end{array}\right.
\end{equation}
For $k=n$ we have the differential operators $\Delta_{R,nj}:C^{\infty}(M,\Lambda_{j}^{n})\rightarrow C^{\infty}(M,\Lambda^{n}_{j})$, $j=1,2$, 
given by the formulas, 
\begin{equation}
    \Delta_{R,n1}= (d_{R,n-1}d^{*}_{R,n})^{2}+D_{R,n}^{*}D_{R,n}, \quad   \Delta_{R,n2}=D_{R,n}D_{R,n}^{*}+  (d^{*}_{R,n+1}d_{R,n}).
         \label{eq:contact-Laplacian2}
\end{equation}

Observe that $\Delta_{R,k}$, $k\neq n$, is a differential operator order $2$, whereas $\Delta_{Rn1}$ and $\Delta_{Rn2}$ are differential operators of 
order $4$. Moreover, Rumin~\cite{Ru:FDVC} proved that in every degree the contact Laplacian is maximal hypoelliptic. 
In fact, in every degree the contact Laplacian has an invertible principal symbol, hence admits a parametrix in the Heisenberg calculus 
(see~\cite{JK:OKTGSU},~\cite[Sect.~3.5]{Po:MAMS1}).
 
Let $\Pi_{0}(d_{R,k})$ and $\Pi_{0}(D_{R,n})$ be the orthogonal projections onto $\ker d_{R,k}$ and $\ker D_{R,n}$, and let 
$\Delta_{R,k}^{-1}$ and $\Delta_{R,nj}^{-1}$ be the partial inverses of $\Delta_{R,k}$ and $\Delta_{R,nj}$.
Then as in~(\ref{eq:projections-dbarb-varthetab}) we have 
%
 \begin{gather}
     \Pi_{0}(d_{R,k})=\left\{ 
     \begin{array}{ll}
         1-(n-k-1)^{-1}d^{*}_{R,k+1}\Delta_{R,k+1}^{-1}d_{R,k}, & k=0,\ldots,n-2,\\
          1- d^{*}_{R,n}d_{R,n-1}d^{*}_{R,n}\Delta_{R,n1}^{-1}d_{R,n-1}, & k=n-1,\\
          1- (k-n)^{-1}d^{*}_{R,k+1}\Delta_{R,k+1}^{-1}d_{R,k}, & \text{$k=n,\ldots,2n-1$},
     \end{array} \right. \\
     \Pi_{0}(D_{R,n})=1-D_{R,*}\Delta_{R,n2}^{-1}D_{R,n}.
  \end{gather}

As in each degree the principal symbol of the contact Laplacian is invertible, 
the operators $\Delta_{R,k}^{-1}$, $k\neq n$, and $\Delta_{R,nj}^{-1}$, $j=1,2$  are \psivdos\ of order $-2$ and order $-4$ respectively. 
Therefore, the above formulas for $\Pi_{0}(d_{R,k})$ and $\Pi_{0}(D_{R,n})$ show that these projections are zero'th order \psivdos.

\begin{theorem}[\cite{Po:AofM1}]
   1) $\Res \Pi_{0}(d_{R,k})$, $k=1,\ldots,2n-1$, and $\Res \Pi_{0}(D_{R,n})$ are Heisenberg diffeomorphism invariants of 
    $M$, hence their values depend neither on the contact form $\theta$, nor on the almost complex structure $J$.\smallskip 
    
    2)  These noncommutative residues are invariant under deformations of the contact structure.  
\end{theorem}

\section{Vanishing of Hirachi's invariant}\label{sec:vanishing}
In this section we give simple algebro-geometric arguments proving that Hirachi's invariant $L(S)=-\frac{1}{2}\Res S_{b;0,0}$ 
always vanishes on strictly pseudoconvex CR manifolds of dimension $4m+1$. First, we have: 

\begin{lemma}\label{lem:Vanishing.Lemma1}
   (i)  We have $\Res \Pi_{0}(\dbarbpq)=-\Res \Pi_{0}(\overline{\partial}_{b;p,q+2}^{*})$ when the condition $Y(q+1)$ holds everywhere.\smallskip
   
   (ii)  We have $\Res \Pi_{0}(\dbarbpq)=\Res \Pi_{0}(\overline{\partial}_{b;p,q+2})$ when the condition $Y(q+1)$ and $Y(q+2)$ both hold everywhere.\smallskip
   
   (iii) If $M$ is strictly pseudoconvex and $n$ even we have $\Res S_{b;p,0}=-\Res S_{b;p,n}$.
\end{lemma}
\begin{proof}
   Suppose that the  condition $Y(q)$ holds everywhere. Then from~(\ref{eq:projections-dbarb-varthetab}) 
   and the fact that $\Res$ is a trace vanishing on $1$ we get 
   \begin{gather}
       \Res \Pi_{0}(\dbarbpq)=-\Res (\overline{\partial}_{b;p,q+1}^{*}N_{b;p,q+1}\dbarbpq)=-\Res (\dbarbpq\overline{\partial}_{b;p,q+1}^{*}N_{b;p,q+1}),\\
         \Res \Pi_{0}(\overline{\partial}_{b;p,q+2}^{*})=-\Res (\overline{\partial}_{b;p,q+1}N_{b;p,q+1}\overline{\partial}_{b;p,q+2}^{*})= 
        -\Res (\overline{\partial}_{b;p,q+2}^{*}\overline{\partial}_{b;p,q+1}N_{b;p,q+1}).
   \end{gather}
  Therefore, we see that $\Res \Pi_{0}(\dbarbpq) +  \Res \Pi_{0}(\overline{\partial}_{b;p,q+2}^{*}) $ is equal to
    \begin{multline}
    - \Res [(\dbarbpq\overline{\partial}_{b;p,q+1}^{*} + \overline{\partial}_{b;p,q+2}^{*}\overline{\partial}_{b;p,q+1}) N_{b;p,q+1}] 
       = -\Res (\Box_{b;p,q+1}N_{b;p,q+1})\\ =-\Res (1-S_{b;p,q+1}).
  \end{multline}
  Since the condition $Y(q+1)$ holds everywhere the operator $S_{b;p,q+1}$ is smoothing and so we have $\Res (1-S_{b;p,q+1})=0$. 
  It then follows that $\Res \Pi_{0}(\dbarbpq)=-\Res \Pi_{0}(\overline{\partial}_{b;p,q+2}^{*})$. 
  
  Assume now that both conditions $Y(q+1)$ and $Y(q+2)$ hold everywhere. Then $S_{b;p,q+2}$ is smoothing and 
  so from~(\ref{eq:Szego-projections-dbarb-varthetab}) we get 
  \begin{equation}
     \Res \Pi_{0}(\overline{\partial}_{b;p,q+2}) +\Res \Pi_{0}(\overline{\partial}_{b;p,q+2}^{*}) =\Res (1+S_{b;p,q+2})=0.
  \end{equation}
  Hence $  \Res \Pi_{0}(\overline{\partial}_{b;p,q+2})=- \Res \Pi_{0}(\overline{\partial}_{b;p,q+2}^{*})=\Res \Pi_{0}(\dbarbpq)$ as desired. 
  
  Finally, suppose that $M$ is strictly pseudoconvex and that $n$ is even. Then the condition $Y(q)$ holds for $q=1,\ldots,n$. In particular, the 
  conditions $Y(q+1)$ and $Y(q+2)$ hold everywhere simultanously for $q=0,2,\ldots,n-4$. Therefore, by the part (ii) we have 
  $\Res \Pi_{0}(\overline{\partial}_{b;p,0}) =  \Res \Pi_{0}(\overline{\partial}_{b;p,2}) = \ldots =\Res \Pi_{0}(\overline{\partial}_{b;p,n-2})$. Moreover, 
  as the condition $Y(n-1)$ holds,  by the part (i) we have $\Res \Pi_{0}(\overline{\partial}_{b;p,n-2})=-\Res \Pi_{0}(\overline{\partial}_{b;p,n}^{*})$. 
  Since $S_{b;p,0}=  \Pi_{0}(\overline{\partial}_{b;p,0}) $ and $S_{b;p,n}=\Pi_{0}(\overline{\partial}_{b;p,n}^{*})$ it follows that $\Res S_{b;p,0}=-\Res 
  S_{b;p,n}$. 
\end{proof}

Next, assume $M$ strictly pseudoconvex and let $\theta$ be a contact form anihilating $T_{1,0}\oplus T_{0,1}$. Let $X_{0}$ be the Reeb vector 
field of $\theta$ so that $\iota_{X_{0}}\theta=1$ ad $\iota_{X_{0}}d\theta=0$. We endow $T_{\C}M$ with the Levi metric associated to $\theta$, i.e., 
the Hermitian metric $h_{\theta}$ such that: \smallskip 

- The splitting $T_{\C}M=T_{1,0}\oplus T_{0,1}\oplus \C X_{0}$ is orthogonal with respect to $h_{\theta}$;\smallskip

- $h_{\theta}$ commutes with complex conjugation;\smallskip

- $h_{\theta}$ agrees with $L_{\theta}$ on $T_{1,0}$ and we have $h_{\theta}(X_{0},X_{0})=1$.\smallskip

\noindent By duality this defines a Hermitian metric on $\Lambda^{*}T^{*}_{\C}M$, still denoted $h_{\theta}$, and there is a uniquely defined Hodge 
operator $*:\Lambda^{p,q}\rightarrow \Lambda^{n-q,n-p}$ such that 
\begin{equation}
    \alpha \wedge \overline{* \beta}=h_{\theta}(\alpha,\beta) d\theta^{n} \qquad \forall \alpha,\beta \in C^{\infty}(M,\Lambda^{p,q}). 
\end{equation}
The operator $*$ is unitary and satisfies $*^{2}=(-1)^{p+q}$ on $\Lambda^{p,q}$. Moreover, we have
\begin{equation}
   \overline{\partial}_{b;p,q}^{*}=-* \partial_{b;n-q,n-p} *, \qquad  \partial_{b;p,q}^{*}=-* \overline{\partial}_{b;n-q,n-p} *.
     \label{eq:CR.Hodge-duality}
\end{equation}


\begin{lemma}\label{lem:Vanishing.Lemma2}
    Assume that $M$ is strictly pseudoconvex. Then we have $\Res S_{b;0,0}=\Res S_{b;0,n}$. 
\end{lemma}
\begin{proof}
    First, let $\Pi_{0}(\partial_{b;0,0})$ be the orthogonal projection onto the kernel of $\partial_{b;0,0}$. Notice that the operators $\overline{\partial}_{b}$ and 
$\partial_{b}$ in~(\ref{eq:CR.dbarb}) are complex conjugates of each other, i.e., we have 
\begin{equation}
    \overline{\partial}_{b;p,q}\alpha = \overline{\partial_{b;p,q}\overline{\alpha}} \qquad \forall \alpha,\beta \in C^{\infty}(M,\Lambda^{p,q}). 
     \label{eq:CR.complex-conjugates}
\end{equation}
Therefore $\Pi_{0}(\partial_{b;0,0})$ is the 
    complex conjugate of $\Pi_{0}(\overline{\partial}_{b;0,0})=S_{b;0,0}$. As $\overline{S_{b;0,0}}=(S_{b;0,0}^{*})^{t}=S_{b;0,0}^{t}$ and by the results 
    of~\cite{Po:GAFA1} we have $c_{S_{b;0,0}^{t}}(x)=c_{S_{b;0,0}}(x)$, we see that the densities $c_{\Pi_{0}(\partial_{b;0,0})}(x)$ and $c_{S_{b;0,0}}(x)$ 
    agree.
    
    On the other hand, let $\Pi_{0}(\overline{\partial}_{b;n,n}^{*})$ denote the orthogonal projection onto $\ker \overline{\partial}_{b;n,n}^{*}$. 
    Since by~(\ref{eq:CR.Hodge-duality}) we have  $\overline{\partial}_{b;n,n}^{*}=-* \partial_{b;0,0} *$ we see that $\Pi_{0}(\overline{\partial}_{b;n,n}^{*})= 
    *\Pi_{0}(\partial_{b;0,0}) *$. As $*^{2}=1$ on $\Lambda^{n,n}$ we get
    $\tr_{\Lambda^{n,n}}c_{\Pi_{0}(\overline{\partial}_{b;n,n}^{*})}(x)=\tr_{\Lambda^{n,n}}[*c_{\Pi_{0}(\partial_{b;0,0})}(x)*]= 
    c_{\Pi_{0}(\partial_{b;0,0})}(x)$. Therefore, we have  $\tr_{\Lambda^{n,n}}c_{\Pi_{0}(\overline{\partial}_{b;n,n}^{*})}(x)= c_{S_{b;0,0}}(x)$, so 
    that $\Res \Pi_{0}(\overline{\partial}_{b;n,n}^{*})=\Res S_{b;0,0}$.

    Next, let $Z_{1},\ldots,Z_{n}$ be a local orthonormal frame of $T_{1,0}$. Then $\{X_{0},Z_{j},Z_{\bar j}\}$ is an 
    orthonormal frame of $T_{\C}M$. Let $\{\theta, \theta^{j},\theta^{\bar j}\}$ be the dual coframe on 
    $T_{\C}^{*}M$. Let $\zeta=\theta^{1}\wedge \ldots \wedge \theta^{n}$ and let $\tau(\eta)=\zeta\wedge\eta$. Then $\tau$ is a 
    a (locally defined) vector bundle isomorphism from $\Lambda^{0,n}$ onto $\Lambda^{n,n}$ and we have 
    \begin{gather}
        \dbarb (\zeta\wedge\eta)= (\dbarb \zeta)\wedge \eta +(-1)^{n}\zeta \wedge \dbarb \eta,\\
        \dbarb \zeta = \sum_{1\leq j\leq n}(-1)^{j-1} \theta^{1}\wedge \ldots \wedge \theta^{j-1}\wedge (\dbarb \theta^{j})\wedge \theta^{j+1}\wedge \ldots. 
        \wedge \theta^{n}.
         \label{eq:CR.dbarb-zeta}
    \end{gather}
    
    Let $\omega \in C^{\infty}(M, T^{*}M \otimes \End(T_{\C}^{*}M))$ be the connection 1-form of the Tanaka-Webster connection (see~\cite{Ta:DGSSPCM}, 
    \cite{We:PHSRH}). Thus, if we let $\omega_{k}^{j}=h_{\theta}(\omega(Z_{j}), Z_{k})$ then 
    $d\theta^{j}=\theta^{k}\wedge \omega_{k}^{j} \ \bmod \theta \wedge T^{*}M$. Let 
   $\omega^{0,1}$ and $\omega_{k}^{j,0,1}$ 
    be the respective $(0,1)$-components of $\omega$ and $\omega_{k}^{j}$. Then we have $\dbarb \theta^{j}=\theta^{k} \wedge 
    \omega_{k}^{j,0,1}$. Combining this with~(\ref{eq:CR.dbarb-zeta}) then gives 
    \begin{equation}
        \dbarb \zeta = \sum_{1\leq j\leq n} (-1)^{j-1}\theta^{1}\wedge \ldots \wedge \theta^{j-1}\wedge \theta^{j} \wedge 
    \omega_{j}^{j,0,1}\wedge \theta^{j+1}\wedge \ldots  \wedge \theta^{n} = -\Tr \omega^{0,1} \wedge \zeta. 
    \end{equation}
    Hence $\dbarb (\zeta\wedge\eta)= (-1)^{n} \zeta \wedge \dbarb \eta - \Tr \omega^{0,1} \wedge \zeta \wedge \eta$. Thus,
    \begin{equation}
        \tau^{-1} \overline{\partial}_{b;n,q}\tau= D_{b;0,q}, \qquad 
        D_{b;0,q}\eta=(-1)^{n}\overline{\partial}_{b;0,q}\eta-\Tr \omega^{0,1}\wedge \eta.
    \end{equation}
    
   As $\tau$ is a unitary isomorphism we also have $\tau^{-1} 
    \overline{\partial}_{b;n,q}^{*}\tau= D_{b;0,q}^{*}$. Using~(\ref{eq:projections-dbarb-varthetab}) 
    we then deduce 
    that $\tau^{-1} \Pi_{0}(\overline{\partial}_{b;n,n}^{*})\tau$ agrees with the orthogonal projection $\Pi_{0}(D_{b;0,n})$ onto 
    $\ker  D_{b;0,n}$. 
    Therefore, the density $ \tr_{\Lambda^{0,n}}c_{ \Pi_{0}(D_{b;0,n})}(x)$ is equal to
    \begin{equation}
        \tr_{\Lambda^{0,n}}[\tau(x)c_{ 
        \Pi_{0}(\overline{\partial}_{b;n,n}^{*})}\tau(x)^{-1}]= \tr_{\Lambda^{n,n}}c_{ \Pi_{0}(\overline{\partial}_{b;n,n}^{*})}(x). 
    \end{equation}
    Hence $\Res \Pi_{0}( D_{b;0,n})= \Res \Pi_{0}(\overline{\partial}_{b;n,n}^{*})$.

   On the other hand, it also follows from~(\ref{eq:projections-dbarb-varthetab}) 
   that $ \Pi_{0}( D_{b;0,n})$ and $ 
   \Pi_{0}(\overline{\partial}_{b;0,n}^{*})$ have same principal symbol, so by~\cite[Prop.~3.7]{Po:AofM1} their noncommutative residues 
   agree. 
   Hence $\Res \Pi_{0}(\overline{\partial}_{b;n,n}^{*})=\Res \Pi_{0}(\overline{\partial}_{b;0,n}^{*})$. 
   As $S_{b;0,n}=\Pi_{0}(\overline{\partial}_{b;0,n}^{*})$ and we have shown above that $\Res \Pi_{0}(\overline{\partial}_{b;n,n}^{*})=\Res S_{b;0,0}$, we 
   see that $\Res S_{b;0,n}=\Res S_{b;0,0}$. The lemma is thus proved.
\end{proof}

We are now ready to prove:

\begin{proposition}
    The Hirachi invariant vanishes on strictly pseudoconvex CR manifolds of dimension $4m+1$.
\end{proposition}
\begin{proof}
    Let $M$ be strictly pseudoconvex CR manifolds of dimension $4m+1$. 
    By Lemma~\ref{lem:Vanishing.Lemma1} we have $\Res S_{b;0,0}=-\Res S_{b;2m,2m}$ and by Lemma~\ref{lem:Vanishing.Lemma2} we have $\Res 
    S_{b;0,0}=\Res S_{b;2m,2m}$, so $\Res S_{b;0,0}=0$. Hence the result. 
\end{proof}

\end{document}